\documentclass{article}
\usepackage{amsmath,latexsym,amssymb,amsthm,enumerate,amscd}
\theoremstyle{plain}
\newtheorem{theorem}{Theorem}[section]
\newtheorem{proposition}[theorem]{Proposition}
\newtheorem{corollary}[theorem]{Corollary}
\newtheorem{lemma}[theorem]{Lemma}
\theoremstyle{definition}
\newtheorem{definition}[theorem]{Definition}

\newtheorem{remark}[theorem]{Remark}
\newtheorem{example}[theorem]{Example}

\newtheorem*{urem}{Remark}
\newtheorem{thm}{Theorem}

\theoremstyle{plain}

\newtheorem{rem}[thm]{Remark}

\newtheorem*{uthm}{Theorem}

\newtheorem*{ack}{Acknowledgment}

\numberwithin{equation}{section}
\numberwithin{table}{section} 
\newcommand{\Gor}{\rm{Gor}}

\DeclareMathOperator{\Grass}{\rm{Grass}}

\def\M{\mathrm{Mat}}

\def\LA{\mathrm{LA}}

\def\cha{\mathrm{char}\ }

\def\Tor{\mathrm{Tor}}

\def\Hilb{\mathrm{Hilb}}

\def\Soc{\mathrm{Soc}}

\def\Grass{\mathrm{Grass}}

\providecommand{\bysame}{\makebox[3em]{\hrulefill}\thinspace}
\def\cod{\mathrm{cod}\ }
\def\<{\left<}
\def\>{\right>}

\def\G{\mathrm{G}}

\def\ns{\footnotesize \it}

\def\ST{\mathrm{ST}}

\def\max{\mathrm{max}}

\def\Anc{\mathrm{Anc}}
\title{Betti strata of height two ideals}
\author{Anthony Iarrobino\\[.05in]
{\ns Department of Mathematics, Northeastern University, Boston, MA 02115, USA.
}\\[.2in]}

\date{November 27, 2004}

\begin{document}

\maketitle
\begin{abstract} Let $R=k[x,y]$ denote the polynomial ring
in two variables over an infinite field $k$.  We study the Betti strata of the family $\G(H)$ parametrizing
 graded Artinian quotients of
$R=k[x,y]$ having given Hilbert function $H$.  The Betti stratum $\G_\beta (H)$
parametrizes all quotients $A$ of having the graded Betti numbers determined by $H$ and the minimal relation
degrees $\beta$, with $\beta_i=\dim_k\Tor_1^R(I,k)_i$. We recover
that the Betti strata are irreducible, and we calculate their codimension
in the family
$\G(H)$.
\begin{uthm}The codimension of $\G_\beta(H)$ in $\G(H)$ satisfies, letting
$\nu_i=\# \{\text{ generators of degree } i\}$, and $\beta_i=\#\{ \text{
relations of degree } i\},$
\begin{equation*}
\cod \G_\beta (H)=\sum_{i>\mu}(\beta_i)(\nu_i).
\end{equation*}
\end{uthm}\noindent
 Here $\mu=\min\{i\mid H_i<i+1\}$, the initial degree of the ideals. When $k$ is algebraically closed, we
also show that the closure of a Betti stratum is the union
of more special strata, and is Cohen-Macaulay (Theorem
\ref{tauideal}). Our method is to identify $\G_\beta (H)$ as a product of determinantal varieties. Ours is a more direct argument than that of \cite[Theorem
5.63]{IK}, which obtains the codimension of Betti substrata of a postulation stratum
for punctual subschemes of
${{\mathbb P}}^2$; there we relied on a result of M.~Boij determining the codimension of Betti
strata of height three Gorenstein algebras \cite{Bj}. Key
tools throughout include properties of an invariant
$\tau (V)$, the number of generators of the ancestor ideal
$\overline{V}\subset R$, and
previous results concerning the projective variety $\G(H)$ in
\cite{I1}.  We adapt a method of M.~Boij that reduces the
calculation of the codimension of the Betti strata to showing that the most special stratum
has the right codimension. \par
 As applications we determine which Hilbert functions are possible for Artinian quotients of $R$, given the socle type (Theorem \ref{4.6C}), and we also determine
the Hilbert function of the intersection of $t$ general enough level ideals in $R$, each having a specified Hilbert function (Theorem \ref{intersectlevelideals}). 
\end{abstract}
\section{Introduction}
Let $R=k[x,y]$ be the polynomial ring in two variables over an infinite field $k$; for our main results we
assume that $\cha k=0$ or $\cha k$ is large enough -- greater than $j=\max\{i\mid H_i\not=0\}$.
We study the Betti strata $\G_\beta (H)$ of the family $\G(H)$ of standard graded 
Artinian quotients $A=R/I$ of $R$ having given Hilbert function $H$; here $\beta$
 is the sequence specifying the number of minimal relations in each degree 
for the ideal $I$, so $\beta_i=\dim_k\Tor_1^R(I,k)_i$ 
\smallskip\par
 For arbitrary graded ideals in
$R$, the Hilbert function does not determine the graded Betti numbers, in
contrast to the situation for the three special algebras determined by a
vector space $V$ of degree-$j$ forms in $R$ studied in \cite{I4} --- the
algebra
$R/(V)$, the level algebra $LA(V)$ and the ancestor algebra $\Anc(V)$.
Which graded Betti numbers are possible, given the Hilbert function is
well known (see \cite{Cam},\cite[\S 5]{F}, and \cite[Lemma 4.5]{I1}). By
the Hilbert-Burch theorem the Betti stratum
$\G_\beta(H)$ is determined by the Hilbert function
$H=H(R/I)$, together with the sequence $\beta=\beta (I)$ specifying the number
of minimal relations in each degree.  A key invariant we used in studying the special
algebras, and that we use again here is $\tau (V)$, a measure of the complexity of
$V$: it is the number of generators of the ancestor ideal 
$\overline{V}$ (see \cite{I4} and Lemma \ref{taugen} below). \par
We will show that, given $H$, to specify 
$\beta (I)$ is the same as to specify the invariant $\tau(I)$, which fixes $\tau_i=\tau (I_i)$ for each
$i$:  and we
give sharp bounds on these sequences of integers (Proposition
\ref{tauideala}).
 We then show that each Betti stratum $\G_\beta(H)$ is irreducible, and we
determine the codimension of $\G_\beta(H)$ in the projective variety $\G(H)$
parametrizing all graded ideals of
$R$ such that
$H(R/I)=H$. Letting $\beta_i,\nu_i$ denote the number of minimal degree-$i$
 relations and generators,
 respectively of $I$, and $\mu=\mu (H)=\min\{i\mid H_i<i+1\}$ the \emph{order} of $H$ (initial degree
 of the ideals $I$), we show (Theorem \ref{tauideal}\eqref{tauideali},
equation \eqref{codBettib}) that this codimension has the strikingly simple form
\begin{equation*}
\cod \G_\beta (H)=\sum_{i>\mu}(\beta_i)(\nu_i).
\end{equation*}
 We also show that the Zariski closure of a Betti stratum satisfies the \emph{frontier
property},  
\begin{equation*}
\overline{\G_\beta (H)}=\bigcup_{\beta'\ge\beta} \G_{\beta'}(H),
\end{equation*} 
where $\beta '\ge \beta$ if for each $i>\mu (H), \beta '_i\ge \beta_i$; also, the closure is Cohen-Macaulay 
with singularity locus $\bigcup_{\beta'>\beta} \G_{\beta'}(H)$. (Theorem \ref{tauideal}\eqref{tauidealii}, 
\eqref{tauidealiii}). \par
Our method is to identify an open subset of the $\G_\beta (H)$ stratum as the product of determinantal varieties,
where each determinantal variety is the rank $\tau_i$ locus of a certain homomorphism 
$\theta_i$ (see \eqref{monomb}).
Here, we mean by \emph{determinantal variety} a variety defined by suitable minors of a matrix,
each irreducible component of which
 has the
right codimension, namely 
\begin{equation}\label{eqdetcod}
(u-\tau)á(v-\tau)
\end{equation}
 for the rank $\tau$ locus of a $u\times v$ matrix (see \cite{BV,ACGH,W}). \par
As a consequence of our main result, we determine in Section \ref{fixed} the Hilbert functions possible for
algebras of given  socle type when
$r=2$, a result announced in \cite[Theorem 4.6C]{I3}, whose complete
proof was referenced to the paper \cite{I4}, but is instead shown here
(Theorem \ref{4.6C}). We then determine the Hilbert function
of the intersection of a set of general enough level ideals $I(i)\subset R$ each
having given invariants  ($\tau$, order, type) and each having a given
Hilbert function
$H(i)=H(R/I(i))$ (Theorem
\ref{intersectlevelideals}). 
\smallskip \par 
That we can show a codimension formula and frontier property for Betti strata of height
two Cohen-Macaulay algebras,
is not surprising nowadays, as similar results have been
found for the Betti strata of height three graded Gorenstein ideals
\cite{Bj,D}.  
 Also, V.~Kanev and the author using the results of S.~J.~Diesel and
M.~Boij derived analogous codimension formulas and showed a frontier
property for the Betti substrata of the schemes
$\Hilb^{\Sigma(H)}({\mathbb P}^2)$, which are the postulation Hilbert
schemes parametrizing those punctual schemes in
${\mathbb P}^2$ with Hilbert function $\Sigma (H)$ \cite[Theorem
5.63]{IK}. The postulation Hilbert function  $\Sigma (H)$ there is the
sum function  of
$H$ here. \par
There, we used as input M.~Boij's formula for the
codimension of Betti strata for height three Gorenstein ideals of a
given Hilbert function $T$ containing a subsequence $(t,t,t)$, and also
S.~J.~Diesel's proof of a frontier property for those strata
\cite{Bj,D}. We also used the connections between certain families of
codimension three Gorenstein Artinian algebras and punctual schemes in
${\mathbb P}^2$, studied extensively in \cite[Part II]{IK}. The
codimension result here Theorem \ref{tauideal}\eqref{tauideali} could be derived from \cite[Theorem 5.63]{IK}.
The argument here is rather more direct.
Conversely, Theorem~\ref{tauideal} below and the argument of
V. Kanev and the author in \cite[\S 5.5]{IK} could be used to recover the
codimension formula of M.~Boij for the Betti strata of height three,
socle-degree $j$ Gorenstein ideals of Hilbert function $T$, with $(\Delta
T)_{\le j/2}=H$,
\cite{Bj}; however, this recovers the M.~Boij result only for the special case that
$T$ contains a subsequence
$(t,t,t)$, where $t$ is the length $|H|=\sum H_i$.\par The
irreducibility of the Betti strata here is a natural consequence of the
Hilbert-Burch theorem giving a determinantal minimal resolution for
height two Cohen-Macaulay ideals \cite{Bu}. G.~Ellingsrud used a similar approach
\cite[p.427]{Ell} to show the irreducibility of the union of the Betti
strata with
$\beta$ smaller or equal a given $\beta_0$ (so more general than a given
stratum), as part of his proof of the smoothness of the Hilbert scheme
of codimension two ACM ideals. The
frontier property here would also admit of a proof analagous to that of
S.~J.~Diesel for the Betti strata of Gorenstein height three ideals, by
deforming the relation matrix 
\cite{D}.\par However, we have chosen to base our
proof of Theorem
\ref{tauideal}  on the relation matrix for the standard generators found
in
\cite{I1} for height two graded Artinian quotients of $R$.
 More
precisely, we base it on the study of a certain homomorphism $\theta_i$, which for each degree $i,\mu\le i<j$
 has rank $\tau_i=\tau (I_i)$, the dimension of the 
vector space quotient $R_1I_i/yI_i$. We feel that the interpretation
of the Betti strata as the product of the rank $\tau_i$
strata of these homomorphisms illuminates the codimension formulas
\eqref{codBettia},\eqref{codBettib}, and also allows us to apply the theory of determinantal ideals.  A key step in the proof uses
 a method of M.~Boij, that determines the
codimension of each Betti stratum ---  a putative product of
determinantal varieties --- by verifying the codimension of the most
special stratum \cite[Proposition 3.2]{Bj}. Here the most special
stratum satisfies, $i\ge
\mu (H) $ implies $ \tau (I_i)=1$: this is equivalent to $I_i$ being the degree-$i$ homogeneous component of a principal ideal
of degree $H_i$ (that depends on $i$),
and is the case $\beta=\beta_\max (H)$
of \eqref{eqbetamax} and Corollary \ref{betamax}.
\par
In higher dimensions, several authors have studied the postulation of height two
irreducible arithmetic Cohen-Macaulay curves in $\mathbb P^3$, in particular showing that they
have ``decreasing type'' \cite{GrP,MaR,HTV,GM}.
In Remark~\ref{higheremb} we discuss the consequences of our main result for Betti strata of height two ACM
schemes of given $h$-vector $H$ and higher embedding dimension than two. 
\section{Betti strata for graded ideals of $R=k[x,y]$}\label{relatedapply}
 We assume throughout that $R=k[x,y]$
is the polynomial ring in two variables over an arbitrary field $k$, unless $k$ is otherwise specified. We study standard graded Artinian quotients $A=R/I$ having given Hilbert function $H$.
In Section \ref{tausect} we show that, given $H$, to specify a Betti stratum $G_\beta (H)$ is 
equivalent to fixing $\tau (I_i)$ for each relevant degree $i$ (Proposition \ref{tauideala} and Corollary
\ref{tauandbetti}). In Section \ref{tauideals} we study the codimension and Zariski
 closure of the Betti stratum
$\G_\beta (H)$. Our main result is Theorem
\ref{tauideal}.  \par
 We first define the invariant $\tau (V)$ and state several preparatory results. 
\begin{definition}\label{def2.1} For $a>0$ and a vector subspace $V\subset R_j$ we let  
\begin{align}  
R_aV&=\langle fv\mid f\in R_a, v\in V\rangle\subset R_{a+j}\\
R_{-a}V&=\{ f\in R_{j-a}\mid  R_af\subset V\}.
\end{align}
We define $\tau (V)$ by
\begin{equation}\label{etauV}
\tau (V)=\dim_k R_1V-\dim_k V,
\end{equation}
and the \emph{ancestor ideal}
\begin{equation}\label{eanc}
\overline{V}=R_{-j}V+\cdots +R_{-1}V+(V)
\end{equation}
\begin{lemma}\cite[Lemma 2.2]{I4}\label{taugen}
We have 
\begin{align}
\tau (V)&=\dim_k V-\dim_k R_{-1}V\text { and }\label{taubasici}\\
 \tau(V)&= \# \{ \text{generators of } \overline{V}\}.\label{taubasicii}
\end{align}
\end{lemma}
\begin{proof}[Remark on proof] The equation \ref{taubasici} results from the exact sequence
\begin{equation}\label{tauseq}
   0\to R_{-1}V\xrightarrow{\phi} R_1\otimes V\xrightarrow{\theta} 
R_1\cdot V
\to 0,
\end{equation}
where $\phi : f\to y\otimes xf-x\otimes yf$, and $\theta
:\sum_i\ell_i\otimes v_i\to \sum_i \ell_iv_i$, where the $\ell_i$ are
elements of $R_1$ (linear forms). The equation \ref{taubasicii} results from a repeated
application of \eqref{tauseq}.
\end{proof}
\end{definition}\noindent
We recall the characterization of $O$-sequences $H$ --- the sequences possible for
Hilbert functions $H(R/I)$ in two variables. We let $\mu
=\min\{ i\mid H_i\ne i+1\}$.
\begin{lemma} The sequence $H$ is an $O$-sequence for some quotient of $R=k[x,y]$ if and only if
\begin{align}
H&=(1,2,\ldots
,\mu,H_\mu,H_{\mu +1},\ldots ,H_i,\ldots ) \text{ where }\notag\\
&\mu\ge H_\mu\ge H_{\mu +1}\ge \cdots \ge c_H\text{ and }\lim_{i\to
\infty}H_i=c_H\ge 0.\label{Hcond}
\end{align} 
\end{lemma}
\noindent
We let $s(H)=\min\{ i\mid H_i=c_H\}$. In the Artinian case $c_H=0$ and $s(H)=j(H)+1$.
Recall that the variety $\G(H)$ parametrizes graded ideals in $R$ of
Hilbert function $H(R/I)=H$. It is defined as a closed subvariety 
\begin{equation}\label{eGH}
\G(H)\subset \prod_{\mu\le i\le s}\Grass(i+1-H_i, R_i)
\end{equation}
of the product of Grassmanians parametrizing sequences of dimension $(i+1-H_i)$ subspaces $V_i\subset
R_i$: here $\G(H)$ parametrizes those sequences satisfying $R_1V_i\subset
V_{i+1}$ for $\mu\le i\le s-1$.  That this is a closed condition is evident, since for each $i$ it is the intersection of the closed conditions 
$xV_i\subset V_{i+1}$ and $yV_i\subset V_{i+1}$. See \cite{I1,Klp,IY}.\par
When $c=c_H>0$, we denote by $H:c$ the sequence satisfying $(H:c)_i=H_{i+c}-c$ for $i\ge 0$. Then we have the well-known
\begin{lemma}\label{commonfactor}(See \cite[Corollary 2.14]{I4}) Let $I\subset R=k[x,y]$ be an ideal
satisfying
$H(R/I)=H, c=c_H>0 $. Then $\exists f\in R_c$ such that $ I=f\cdot I'$
 where $R/I'$
is an Artinian quotient of $R$ having Hilbert function $H:c$. Also, $\G(H)\cong G(H:c)\times \mathbb P^c$.
\end{lemma}
Evidently, the generators and relations for $I$ are simply related to those for $I'$, and the Betti stratum
$\G_\beta (H)$ of $\G(H)$ is the product of a corresponding Betti stratum for $\G(H')$ by $\mathbb P^c$.
Hence it will 
be sufficient to treat in detail the case $c_H=0$ below.\par
 We will use the
following result, quoted from \cite[Theorem 1.10]{I4}, valid over a field
$k$ of
arbitrary characteristic.  The proof is essentially from
\cite{I1}, but the necessary adjustments for an arbitrary field $k$ are outlined in \cite[Theorem 1.10]{I4}.
Given an $O$-sequence $H$ we set $e_i=e_i(H)=H_{i-1}-H_{i}$. For
an integer $n\in \mathbb Z$ we set $n^+=\max\{ 0,n\}$.
\begin{theorem}\label{basic}\cite[Theorems 2.9,2.12,3.13,4.3, Proposition
4.4, Equation 4.7]{I1}, \cite[Theorem 1.10]{I4} Let
$r=2$, and for \eqref{basici} let the field $k$ be
algebraically closed. Let 
$H$ be an
$O$-sequence that is eventually constant, so $H$ is a sequence satisfying
\eqref{Hcond}, let $c=c_H$ and let $H_s=c_H, H_{s-1}\not=c_H$.
\begin{enumerate}[i.] 
\item\label{basici} Then $\G(H)$ is a smooth projective variety of
dimension $c+\sum_{i\ge \mu} (e_i+1)(e_{i+1})$. $\G(H)$ has a finite
cover by opens in an affine space of this dimension. If $\cha k=0$ or 
$\cha k\ge s$ then $\G(H)$ has a finite cover by opens that are affine spaces.
\item\label{basicii} \cite[Theorem 4.3]{I1} The number of
generators
$\nu (I)$ of a graded ideal
$I$ for which
$H(R/I)=H$, satisfies $\nu (I)\ge \nu (H)=1+e_\mu +\sum_{i\ge
\mu}(e_{i+1}-e_i)^+$.
\item\label{basiciii}\cite[Proposition 4.4]{I1} Assume
that
$k$ is an infinite field. 
The family $U(H)$ of graded ideals
$I$ satisfying
$H(R/I)=H$ and having the minimum possible number $\nu (H)$ of generators given by
equality in
\eqref{basicii} form an open subscheme of $\G(H)$ having the dimension
specified in \eqref{basici}. Also $U(H)$ is dense in $\G(H)$ when $k$ is
algebraically closed.
\end{enumerate}
\end{theorem}
See also Theorem \ref{stdbasis} below for the affine space parameters mentioned in 
Theorem \ref{basic}\eqref{basici}.

\subsection{The $\tau$ invariant and the Betti strata}\label{tausect}

We now specify the relationships between the invariants $\beta,\nu$, and $\tau$ for quotient algebras $A$ of $k[x,y]$.
\begin{definition}\label{defEH} Let $H$ be an $O$-sequence as in
\eqref{Hcond}. We define
$E(H)=(e_\mu,\ldots ,e_i,
\ldots )
$ to be the first difference sequence $e_i=H_{i-1}-H_i$. 
Recall that $s=s(H)=$ $\min\{i\mid H_i=c_H\}$. We denote by $\beta(I)$ the sequence of
integers
\begin{equation}\label{eqbeta}
\beta(I)=(\beta_{\mu+1}(I),\ldots , \beta_i(I),\ldots ,\beta_{s}(I)),
\end{equation}
where $\beta_i(I)$ is the minimal number of relations of $I$ in
degree $i$ (so $\beta_i$ is uniquely determined by $I$). We denote by
$\nu(I)=(\nu_\mu,\ldots ,\nu_s)$ the corresponding sequence where $\nu_i$
is the number of degree-$i$ generators of
$I$, namely $\nu_i=\Tor_0(I_i,k)$, and by
$\tau (I)$ the sequence (see \eqref{taubasici})
\begin{equation}
\tau(I)=(\tau (I_\mu),\tau (I_{\mu
+1}),\ldots ,\tau (I_s)=1)=(\tau_\mu,\tau_{\mu+1},\ldots ,\tau_s ).
\end{equation} 
\end{definition}
We will need the well known fact for the socle $0:M\subset A$ of a graded
quotient
$A=R/I, R=k[x,y]$
\begin{equation}\label{Socandrel}
\dim_k \Soc
(A)_i=\beta (I)_{i+2}.
\end{equation}
The following result gives bounds, in terms of $E(H)=\Delta H$, on the sequences $\nu (I),\beta(I),\tau(I)$. The main result
Theorem \ref{tauideal} below shows that each sequence satisfying the bounds actually occurs. The characterization of $\beta(I)$ 
given $H$ was well known. For example, \eqref{enogensI} follows from
\cite[Lemma 4.5]{I1}, and $\nu$ can be read from $\beta$ and $H$ using \eqref{etauequal}. The characterization of $\beta$ is also a consequence of R. Fr\"{o}berg's
result that every set of Betti numbers for a CM height two graded ideal in $R$
is attained by a monomial ideal \cite[\S 5]{F}. Indeed, the possible
Betti sequences for height two CM ideals were explicitly determined by G. Campanella
in \cite{Cam}.\par 
\begin{proposition}\label{tauideala}
Let $k$ be an artibrary field, and $H=H(R/I)$ where $I$ is a graded ideal of $R=k[x,y]$ be a sequence
satisfying \eqref{Hcond}, and
let
$\mu (H)=\min
\{ i\mid H_i<i+1\}$. Let
$i\ge
\mu (H)$. Then the minimal number of generators $\nu_{i+1}(I)$ of $I$ having
degree $i+1$ satisfies 
\begin{equation}\label{enogensI}
e_{i+1}\ge\nu_{i+1}(I)\ge (e_{i+1}-e_i)^+;
\end{equation}
and $\nu_\mu (I)=e_\mu+1$. 
For $i\ge \mu$ the minimal number of degree $i+1$ relations
$\beta_{i+1}(I)$ satisfies 
\begin{equation}\label{enorelsI}
e_i\ge\beta_{i+1}(I)\,\,\ge (e_i-e_{i+1})^+,
\end{equation} 
and $\beta_i=0$ for $i\le \mu$.
The $\tau$ invariant is zero for $i<\mu$;  for $i\ge \mu$ it satisfies
\begin{align}\label{etauequal}	
\tau
(I_i)&=e_{i+1}(H)+1-\nu_{i+1}(I)=e_i(H)+1-\beta_{i+1}(I)\text { and
}\\\label{etauequalb}
\tau (I_i)&=\tau (
I_{i-1})+\nu_i(I)-\beta_{i+1}(I),
\end{align}
and also for $i\ge \mu$ the inequality
\begin{equation}\label{etauinequal}
1\le\tau (I_i)\le 1+ \min\{ e_i(H), e_{i+1}(H)\}.
\end{equation}
For an open dense subset of $\G(H)$ the 
invariants
$\nu (I), \beta (I),$ and $\tau (I)$, respectively, satisfy equality on the right sides
of
\eqref{enogensI},\eqref{enorelsI} and \eqref{etauinequal}, respectively, 
\end{proposition}
\begin{proof}
We have evidently $\dim_k I_{i+1}=\dim_k R_1I_i+\nu_{i+1}(I)$, and also as a
consequence of \eqref{Socandrel} that $\dim_k R_{-1}I=\dim_k
I_{i-1}+\beta_{i+1}(I)$. Thus by the definition (see
\eqref{etauV}) of
$\tau$ and Lemma
\ref{taugen} we have for $i\ge \mu$,\begin{align*}
\tau (I_i)&=\dim_k R_1I_i-\dim_k I_i= \dim_k I_{i+1	}-\nu_{i+1}(I)-\dim_k
I_i=e_{i+1}+1-\nu_{i+1}(I)
\text { and }\\
\tau (I_i)&=\dim_k I_i-\dim_k R_{-1}I_i= \dim_k I_i-(\dim_k I_{i-1}+\beta_{i+1}(I))
=e_i+1-\beta_{i+1}(I).
\end{align*}
Here the second equation follows from the first, and the usual calculation
that the second difference 
\begin{equation}\label{Hgenrel}
\Delta^2 (H)_{i+1}=e_i-e_{i+1}=\beta_{i+1}(I)-\nu_{i+1}(I).
\end{equation}
This
shows \eqref{etauequal} and implies all the equations \eqref{enogensI} to
\eqref{etauinequal} for
$i\ge
\mu (H)$. For example, \eqref{Hgenrel} directly implies the right sides of
\eqref{enogensI} and \eqref{enorelsI}, while $\tau (I_i)\ge 1$ from \eqref{taubasicii}
together with \eqref{etauequal} imply the left sides of
\eqref{enogensI} and \eqref{enorelsI}; and \eqref{etauequal}
directly implies the right side of \eqref{etauinequal}.  By
Theorem~\ref{basic}\eqref{basiciii}, an open dense subset of ideals $I\in
\G(H)$ consists of ideals having the minimum possible number of generators, given by
$\nu(H)=1+e_\mu+\sum_{i\ge \mu} (e_{i+1}-e_i)^+$. This subset consists of ideals
 satisfying equality on the right of each of the inequalities \eqref{enogensI}, \eqref{enorelsI},\eqref{etauinequal} for $\nu (I), \beta (I)$ and $\tau (I)$, respectively. 
\end{proof}
\begin{corollary}\label{tauandbetti}
Let $I$ be an arbitrary graded ideal
of
$R=k[x,y]$. Then the pair $(H(R/I),\tau (I))$ are equivalent to the pair
$(\nu  (I), \beta (I))$ of sequences
giving the
graded Betti numbers of
$I$. Given $H=H(R/I)$ each of the sequences $\nu(I), \beta(I), \tau(I)$
determines the other
two. 
\end{corollary}
\begin{proof}
We have $\tau_\mu (I)=e_\mu (H)+1=\nu_\mu (I)$, and $\beta_\mu = 0$. Then
\eqref{etauequal} shows that the pair $(H(R/I),\tau (I))$ determines both
the generator and relation degrees, and \eqref{etauequalb} shows that the
latter pair determine the former.
\end{proof}\par
We define the sequence $\beta_\max (H)$ --- the most special stratum --- by equality in the upper bound 
(left side) of \eqref{enorelsI}
\begin{equation}\label{eqbetamax}
(\beta_\max)_{i}=\begin{cases}&e_{i-1}  \text { if } \mu<i\le s(H)\\
&0 \text { otherwise. }
\end{cases}
\end{equation}
The corresponding sequence $\tau_\max$ satisfies $\tau_i=1$ for $i\ge \mu$, and the sequence 
$\nu_\max$ is given by equality on the left of
\eqref{enogensI}. The following result shows that $\beta_\max$ actually occurs.
\begin{corollary}\label{betamax}
Let $H$ be a sequence satisfying \eqref{Hcond}. The subscheme
$\G_{\beta_\max}(H)$ of $\G(H)$ is isomorphic to 
$( \prod_{i=\mu +1}^s{\mathbb P}^{e_i})\cdot {\mathbb P}^{c_H}$,
and has dimension
$H_\mu$.
\end{corollary}
\begin{proof}
When $\beta=\beta_\max$,
$\tau(I)=\tau_{\min}=(1,1,\ldots ,1)$. By Lemma~\ref{taugen} equation \eqref{taubasicii} we have
$\overline{I_i}=(f_i)$ of degree $H_i$, and it is easy to see that
we have $f_s\mid f_{s-1}\ldots \mid  f_\mu$, whence each
such ideal $I$ determines a point $(f_\mu/f_{\mu-1},\ldots ,
f_i/f_{i+1},\dots , f_{s-1}/f_s; f_s)$ of
$(\prod_{i=\mu +1}^s {\mathbb P}^{e_i})\cdot {\mathbb P}^{c_H}$; and conversely.
\end{proof}
\begin{definition}\label{defbettistratum}
Let $H$ be fixed satisfying \eqref{Hcond}, and let $\beta$ be a sequence as in \eqref{eqbeta} satisfying \eqref{enorelsI}. We denote by
$\G_\beta(H)$ the reduced subscheme of $\G(H)$ parametrizing all ideals $I\in R$ such that 
$H(R/I)=H$ and $\beta(I)=\beta$. \par
We will say for two sequences $\beta ',\beta $ satisfying \eqref{enorelsI}
that 
\begin{equation}\label{ebetaineq}
\beta'\ge\beta \text {  if and only if  } \beta_i'\ge \beta_i \,\,\forall i.
\end{equation}
\end{definition}\noindent
Evidently, $\G_\beta (H)$ is a locally closed subscheme of $\G(H)$, as, given $H$, fixing each $\tau_i$ is
to fix the vector space dimension of $R_1I_i$, or, equivalently, the rank of the image of $R_1\otimes I_i\to I_{i+1}$
under the multiplication map.
\subsection{Codimension and closure of the Betti strata}\label{tauideals}
We now state and prove our main result about the Betti strata $\G_\beta(H)$
of
$\G(H)$, the projective variety parametrizing graded quotients of
$R=k[x,y]$ having Hilbert function $H$ (Theorem \ref{tauideal}).
We adapt a result of M.~Boij \cite{Bj}. It follows from equation \eqref{eqbetamax} and
Corollary~\ref{betamax} that for $
\mu(H) \le i< s(H)$ we have $ (\beta_{\max})_{i+1}=e_{i}$. For any ideal
$I\in \G(H)$ we have a resolution
\begin{equation}\label{eqstdburch}
0\to \oplus_{i\ge \mu (H)}R(-i-1)^{e_i}\to R(-\mu )\oplus_{i\ge
\mu(H)}R(-i)^{e_i}\to R\to R/I\to 0,
\end{equation}
that is minimal only when $\beta (I)=\beta_\max$. Denote by $X(H)$ the
open subfamily of all degree zero homomorphisms $\iota$
\begin{equation}\label{eqdefxh}
 \oplus_{i\ge \mu(H)} R(-i-1)^{e_i}\xrightarrow{\iota}
R(-\mu )\oplus_{i\ge
\mu(H)}R(-i)^{e_i}  
\end{equation} 
that have a rank one cokernel.  Such a map $\iota$ corresponding to a point $p\in X(H)$ is
 from a rank $\mu (H)=\sum_{i\ge \mu (H) }e_i$ free module on the left of \eqref{eqdefxh} to one of rank $\mu (H)+1$, so is
injective. By the Hilbert-Burch theorem, the sequence \eqref{eqdefxh} determines an ideal $I=\pi (p)$ of $R$ as in \eqref{eqstdburch},
 defined
by the maximal minors of $\iota$. Thus, we have a morphism $\pi :X(H)\to \G(H)$ taking the point $p\in X(H)$ 
to the ideal $I=\pi (p)$. We denote by $X_\beta (H)$ the
inverse image of $\overline{G_\beta (H)}$. Then, following M.~Boij, we have
a commutative fiber diagram
\begin{equation}
\begin{CD}
 X_\beta (H)  @>>> X(H)\\
   @VVV @VVV                  \\
\overline{G_\beta (H)} @>>> \G(H)
\end{CD}
\end{equation}
where the vertical morphisms are dominant and the horizontal morphisms
are injective. We can replace all the schemes by the associated reduced schemes.
\begin{lemma}\label{boij3.2} (after \cite[Proposition 3.2]{Bj}). For
$\beta$ satisfying \eqref{enorelsI} we have that the codimensions of
$\G_\beta (H)$ in $\G(H)$ and of $X_\beta (H)$ in $X(H)$ satisfy
\begin{equation}\label{codbstrata} \cod \G_{\beta_\max}(H) -\cod
X_{\beta_\max}(H)\ge
\cod \G_\beta (H)-\cod X_\beta (H)\ge 0.
\end{equation}
\end{lemma}
\begin{proof} Since $\pi : X_\beta (H) \to \overline{ \G_\beta (H)}$ is a
dominant morphism of integral schemes the dimension of its general
fibre is no greater than that of its special fibres; the same is true for
$\pi : X(H)\to \G(H)$, hence we have
\begin{equation}
\dim X(H)-\dim \G(H)\le \dim X_\beta
(H)-\dim \G_\beta (H) \le
\dim X_{\beta_\max} (H)-\dim \G_{\beta_\max} (H),
\end{equation}
which is \eqref{codbstrata}.
\end{proof}\par
We now briefly describe $X_\beta (H)$. The standard matrix $\M (H)_p$ for a point $p\in X(H)$, contains an $e_{i+1}\times e_{i}$ submatrix
 $\M _i(H)_p$, having
degree-zero entries, arising from the homomorphism from $R(-i-1)^{e_i}$ to $R(-i-1)^{e_{i+1}}$ in the definition
 \eqref{eqdefxh} of $X(H)$.
\begin{lemma}\label{gencod} Assume that $\beta$ satifies \eqref{enorelsI}. The open subvariety $\pi^{-1}(G_\beta (H))\subset X_\beta(H)$ is the locus of points 
 $p\in X(H)$ for which each
submatrix $\M_i(H)_p$ has rank $e_{i}-\beta_{i+1}=\tau_i -1$, imposing a condition of codimension $\beta_{i+1}\cdot \nu_{i+1}$
 for each such block,
$  \mu (H)\le i\le j$. The locus
 $X_\beta (H)$ has codimension $\sum_{i>\mu ( H)}(\beta_i)(\nu_i)$ in $X(H)$. 
\end{lemma}
\begin{proof} That $\M_i(H)_p$ has rank $e_i-\beta_{i+1}$
leaves $\beta_{i+1}$ uncancelled relations of degree $i+1$ in a minimal generating set for the ideal $\pi (p)$, which is 
determined by the maximal minors of $\M (H)_p$. For a generic point 
$p_{gen} \in X(H)$, the matrix $\M_i(H)_{p_{gen}}$ is generic.
That $\M_i(H)_p$ has rank $e_i-\beta_{i+1}$ is thus a determinantal condition
having the expected codimension 
\begin{equation*}
(e_i-(e_i-\beta_{i+1}))\times (e_{i+1}-(e_i-\beta_{i+1}))=(\beta_{i+1})(\nu_{i+1}),
\end{equation*}
where the equality is by equation \eqref{etauequal}.
\end{proof}\par
 We now begin an explicit description of $G_\beta (H)$ as product of determinantal varieties 
(Lemmas~\ref{ranktau}, \ref{lemgenericmatrix}, Theorem \ref{tauideal}).
We first define an initial monomial ideal $M(H)$ and an open subset $U_B$ of $\G(H)$, the ideals having
\emph{normal pattern} or cobasis (see \cite{I1}). Recall that $s(H)=\min\{ i\mid H_i=c_H\} $.
Recall that given a partition $P$, the dual partition $P^\ast$ is that obtained by exchanging the rows
and columns of the Ferrer's graph of $P$. We use this to define the ``alignment character''
of $H$ below.
\begin{definition}[Normal pattern]\label{defnormal}
Given
$H$ satisfying \eqref{Hcond} we define a \emph{normal pattern} or \emph{standard cobasis}:
\begin{equation}\label{normalpattern}
B=B_H=\oplus_{u=o}^s B_u,\quad B_u=\langle
y^{H_u-1}x^{u+1-H_u},\ldots ,x^u\rangle.
\end{equation}
We say that the ideal $I\in \G(H)$ has \emph{normal pattern} in the coordinates
$x,y$ for $R$ and write $I\in U_B, B=B_H$ if $I$ satisfies
\begin{equation}\label{UB}
 I_u\oplus B_u=R_u\,\,\forall u\mid
\mu\le u\le s(H) .
\end{equation}
We denote by $M(H)$ the unique monomial ideal
$M(H)$ having normal pattern in the basis $H$. It satisfies $M(H)=\oplus M(H)_u$, where
\begin{equation}\label{emonomidH}
   \quad M(H)_u=\langle
y^u,y^{u-1}x,\ldots y^{H_u}x^{u-H_u}\rangle
\end{equation}
\end{definition}
\smallskip\par\noindent
{\it Preparation for standard generators}. We will
need the standard generators for a graded ideal
$I\subset R$ having \emph{normal pattern} or cobasis $B$, and a variation of the coordinates for $U_B$
 used in the
proofs 
\cite[Theorem 4.3 and Proposition 4.4]{I1}. In \cite[Proposition
3.2]{I1} it is shown that given $\cha k=0$ or $\cha k\ge s(H)$ then $\G(H)$ has a finite
cover by affine open subvarieties, each isomorphic to $U_B, B=B_H$. We
will now assume that $I\in U_B$. Our standard generating set $F=F(I)$
for
$I$ has $e_i$ elements of degree
$i$ for
$i>\mu$, and
$1+e_\mu$ elements of degree $\mu$, and it is not usually a minimal
generating set. To present $F$,  following \cite[Proposition 2.5, Lemmas
2.7,2,8]{I1}, we first define a sequence of integers
$K=(k_0>\ldots > k_\mu)$ by adding $(\mu-1,\ldots ,1)$ to the dual
partition
$(h^\ast_1\ge \ldots )$ to
$(H_{\mu-1}=\mu,H_\mu,\ldots ,H_{s-1})$: thus, 
\begin{equation}\label{eqalign}
k_0=h^\ast_1+\mu -1,\ldots , k_i=h^\ast_{i+1}+\mu-(i+1),\ldots .
\end{equation}
The sequence $K$ specifies
the lengths of the horizontal rows of the bar graph of $H$ (see Figure \ref{fig1} and Example
\ref{ex1}), and has sometimes been called the \emph{alignment character} of $H$: see \cite{GPS,Hari}, and also
\cite[p. 122--124]{IK} for a discussion and further references.  For $0\le i\le \mu$,
 the monomial $\upsilon_i  =y^ix^{k_i}$ is the unique monomial
of the ideal $M_H$ having $y$-degree $i$ and lowest total degree. 
\begin{definition}[Standard generators]
For each $I\in U_B$ we let $F=F(I)=(f_\mu ,\ldots ,f_0)$ be the unique homogeneous elements of $I$
such that 
\begin{equation} f_i=y^ix^{k_i}+g_i \text { where }  g_i\subset B_{i+k_i}.
\end{equation}
We define $a_{i,u}\in k$ for $0\le u\le i\le \mu$ by writing
\begin{equation}F: \quad f_i=y^ix^{k_i}+\sum_{\deg F_u>\deg F_i}
a_{i,u}(f_u:x^{k_u-k_i-(i-u)}).
\end{equation}
Note that
$F_\upsilon =(f_{H_{\upsilon -1}-1},\ldots ,f_{H_\upsilon })$, for $s(H)>\upsilon >\mu$, and 
$F_\mu=(f_\mu,\ldots ,f_{H_\mu})$.
comprise the
elements $F_\upsilon =F\cap R_\upsilon$ of $F$ having degree
$\upsilon$ for each
$\upsilon\ge \mu$. We have $\#\{ F_\upsilon\} =e_\upsilon$ when $\upsilon>\mu$ but $\# \{ F_\mu\}=e_\mu+1$,
 and $\#\{ F\}=\mu+1$. \par  We define for $\upsilon \ge\mu$ the set $A_\upsilon$ of coefficients,
in two equivalent ways
\begin{enumerate}[i.]
\item $A_\upsilon =\{ a_{i,u}$ such that $H_\upsilon\ge i\ge H_{\upsilon +1}$ 
and 
$H_{\upsilon +1}-1\ge u\ge H_{\upsilon +2}$\};
\item $A_\upsilon =\{a_{i,u}$ such that $f_i\in F_\upsilon$, or is the next 
generator $f_{H_{\upsilon}}$ of a lower degree, and $f_u
\in F_{\upsilon +1}$\}. We let $A=\bigcup_{\mu\le \upsilon} A_\upsilon$.
\end{enumerate}
\end{definition}
These may also be described as the coefficients of $f_u$ that land outside the pattern $B_H$ in
$x^{w_u}f_u, w_u=k_{u-1}-k_u$ (see \cite{I1} and also \cite[\S 3B]{IY} for a discussion).\par
Evidently $\#\{ A_\upsilon\} =(e_\upsilon +1)e_{\upsilon +1}$.
 Our standard generators differ
from those in \cite[Proposition 2.7]{I1} by using the forms
$f_u:x^a$ here in place of the monomials $y^ux^{i-u}$ used there to write
$f_i$. The chief advantage is to simplify the matrices $N_i$ constructed
below, but there is no change in determining which coefficients are
parameters for
$U_B$.  According to \cite[Lemma 2.7, Theorem
2.9]{I1} these parameters are exactly the set $A$ of coefficients defined above. Thus we have,
\begin{theorem}\label{stdbasis} \cite[Proposition 2.7, Theorem 2.9]{I1} Let $k$ be infinite. 
Fix an $O$-sequence $H$ satisfying \eqref{Hcond} and $c_H=0$, and let $B=B_H$,
 the normal pattern of  \eqref{normalpattern}.
\begin{enumerate}[i.]
\item Each ideal $I\in U_B$ is generated by $F(I)$. 
\item The subfamily $U_B\subset \G(H)$ is isomorphic to an affine space, with coordinates given by the 
coefficients $A$.
\end{enumerate}
\end{theorem}
\begin{proof} That $F(I)$ generates $I$ is evident. Although the coefficients used for $U_B$ differ slightly
from those in \cite[Proposition 2.7, Theorem 2.9]{I1} by using the forms
$f_u:x^a$ here in place of the monomials $y^ux^{i-u}$ used there to write
$f_i$ and hence define $a_{i,u}$; there is no essential change in the result, and in particular, there is
 no change in the definition of the
coefficients used as parameters (those in $A$). Furthermore, each coefficient $a_{i,u}$
 that is not in $A$
may be written as an explicit
polynomial in those parameters in $A$ corresponding to monomials of equal or higher degrees to $\deg f_i$.
 Evidently  $\# \{A\}=
\sum_{i\ge
\mu}(e_i+1)e_{i+1}$, which is the dimension
of $\G(H)$ given in Theorem \ref{basic} above. 
\end{proof}\par
We now prepare for our main result, Theorem \ref{tauideal}, by defining key homomorphisms $\theta_i$.
 Using these we will show that the Betti strata $\G_\beta (H)$ are products of determinantal varieties and have an expected dimension. We assume that $H$ 
has been chosen satisfying \eqref{Hcond} with $c_H=0$, that $I\in U_B$, and that $i$ satisfies $\mu(H)\le i< s(H)$.
\par\smallskip\noindent
{\it The homomorphism $\theta_i$ and the vector space $\langle R_1I_i/yI_i\rangle$}. We have evidently
from the Definition \ref{def2.1} of $\tau$ 
\begin{equation*}
\tau
(I_i)=\dim_k R_1I_i-\dim_k I_i=\dim_k \langle R_1I_i/yI_i\rangle .
\end{equation*}
For now we fix $i$ satisfying  $\mu(H)\le i\le s(H)$ and let $t=H_{i-1}, t'=H_i$. We denote by $F'_i$ the $k$-vector space spanned by
 $(x^{i-\deg f_t}f_t,F_i)$ and thus by
$F'_{i+1}$ the span of $(x^{i+1-\deg f_{t'}}f_{t'},F_{i+1})$. Suppose $I\in U_B$ and consider for $\mu(H)\le i< s(H)$ the homomorphism
\begin{equation}\label{monomb}
\theta_i : F'_i \rightarrow F'_{i+1} ,
\end{equation}
where $\theta_i$ is the composite $\theta_i=\eta_i\circ \theta'_i\circ m_x$ of the homomorphisms
\begin{equation}\label{monoma}
F'_i\xrightarrow{m_x} \langle R_1I_i\rangle
\xrightarrow{\theta '_i}
\langle R_1I_i/yI_i\rangle\xrightarrow{\eta_i} F'_{i+1} ,
\end{equation}
 where $m_x$ denotes multiplication by $x$, the homomorphism $\theta_i'$ is
simply reduction mod $yI_i$, and the inclusion $\eta_i$ arises from the inverse of the
isomorphism
\begin{equation}\label{monomc}
F'_{i+1}\to  \langle F'_{i+1}+yI_i\rangle/yI_i \to I_{i+1}/yI_I;
\end{equation}
The composite map of \eqref{monomc} is an isomorphism, 
 since $I_{i+1}\cap B_{i+1}=0$, and the normal form of $I$ assures that 
$I_{i+1}=yI_i\oplus F'_{i+1}$. 
\begin{lemma}\label{ranktau}
 The
homomorphism $\theta_i$ of \eqref{monomb} is from a vector space of dimension
$e_i+1$ to another of dimension $e_{i+1}+1$; and for $I\in U_{B(H)}$, we have
\begin{equation}\label{tauilocus}
 \text { rank }{\theta_i}=\tau (I_i).
\end{equation}
\end{lemma}
\begin{proof} Since $m_x$ and $\eta_i$
are monomorphisms \eqref{tauilocus}  follows.
\end{proof}\par
We denote by $N_i$ the 
$(e_{i+1}+1)\times (e_i+1)$ matrix for the homomorphism $\theta_i$, where we list both domain and range in 
descending order of the index of the standard generators involved, and where we write the entries in terms of the
parameters $A$ (see Examples~\ref{ex1},\ref{ex2}). We denote by $A'_{u,v}$ the subset of the parameters
involving pairs $f_{u'},f_{v'}$ with $u'<u$, or $u'=u$ and $v'<v$, that is, all those parameters
``earlier than'' the coefficient of $f_u$ on $f_v:x^{k_v-k_u+u-v}$. 
\begin{lemma}\label{lemgenericmatrix}
The entries of $N_i$ are elements of the polynomial ring $k[A]$. The entry $N_i(v,u)$, corresponding to the component of
$\theta_i(f_u)$ involving $f_v$, satisfies
\begin{equation}\label{eqdeterminant}
N_i(v,u)=a_{u,v} \mod k[ A'_{u,v}],
\end{equation}
 for all pairs $u,v$ with $f_v\in
F'_{i+1}$ and $v\ne t'$. Also $N_i$ has last  $f_{t'}$ column $(1,0,\ldots )^T$. After column
reduction (change of basis in the domain of $\theta_i$), the resulting matrix $N'_i$ has
first row $(0,\ldots ,0,1)$.
\begin{proof} It is easy to see inductively, beginning with $f_0$, that for an ideal $I\in U_B$,
 the coefficients of the standard generators $f_i$ are elements of the polynomial ring $k[A]$. From the definition of
the homomorphism $\theta_i$, it follows that the entries of $N_i$ are polynomials
 in the coefficients of the $f_i$ involved, hence are also elements of the polynomial
 ring $k[A]$. Equation \eqref{eqdeterminant} is straightforward, and follows from a
 closer examination of the coefficients of the standard generators involved in calculating the
matrix entry of $N_i$ corresponding to the 
$\theta_i (f_u)$ component on an $f_v, v<u$: this coefficient only involves either ``earlier''
  $f_{u'},u'<u$, or the same $f_u$ but earlier $f_{v'},v'<v$. However, $a_{u,v}$ does occur, as one
evaluates $x^{k_v-k_u+u-v}f_u$ on $f_v$. We have $\theta (f_{t'})=xf_{t'}$, hence the last
column of $N_i$ is $(1,0,\ldots ,0)^T$, as claimed.
\end{proof}
\end{lemma}
\begin{urem} According to Lemma \ref{lemgenericmatrix}, after suitable column reduction (change of basis in the domain of
$\theta_i$)
 the matrices $N_i$ for $\theta_i$ become matrices $N'_i$, which, after removal of the first row and last column
 are generic matrices. Thus, Lemma \ref{lemgenericmatrix} serves to
identify the product of the rank $\tau_i$ loci of the maps $\theta_i$  --- which are the varieties $G_\beta(H)$ --- over a coordinate patch
 (those ideals having
normal cobasis B) as a product of explicit
determinantal varieties. 
 The genericity of the matrices $N'_i$ underlies Theorem~\ref{tauideal} below, and we could have used Lemma \ref{lemgenericmatrix} as the
key to its proof. The proof we
give below of Theorem \ref{tauideal} uses instead M.~Boij's Lemma~\ref{boij3.2} to establish that the varieties 
$\overline{G_\beta (H)}$ are the product of determinantal varieties. 
\end{urem}
 Recall that a 
determinantal variety $X$ is one defined by the rank $\tau$ locus of a linear map or matrix, such that each irreducible component of $X$ 
has the right codimension (see \eqref{eqdetcod}). Recall that we denote by $\G_\beta(H)$ the Betti stratum 
of $\G(H)$ (see Definition \ref{defbettistratum}).
\begin{theorem}\label{tauideal}
Let $H=(1,2,\ldots ,0)$ be an $O$-sequence for an Artinian quotient of $R$: so $H$ satisfies 
\eqref{Hcond} with $c_H=0$. Let $\beta$
denote a sequence of integers (potential relation degrees) satisfying the condition 
\eqref{enorelsI}, and let ${\mathbb \tau}$ and $\nu$ be the sequences defined from $\beta$ by \eqref{etauequal} and 
\eqref{Hgenrel}, respectively.
Let $k$ be an infinite field satisfying either
$\cha k = 0$ or
$\cha k=p\ge s(H)$.  
\begin{enumerate}[i.]
\item\label{tauideali} Then $\G_\beta (H)$ is irreducible. The intersection
$\G_\beta (H)\cap U_B$ is open dense in $\G_\beta (H)$, and is the
product of the rank
$\tau_i$ loci of the homomorphisms $\theta_i$ of \eqref{monomb}, 
determined by the $(e_{i+1}+1)\times (e_i+1)$ matrices $N_i$, $\mu\le i\le s(H)-1$. The codimension in
$\G(H)$ of the Betti stratum $\G_\beta (H)$ satisfies
\begin{align}\label{codBettia}\cod \G_\beta (H)&=\sum_{i\ge\mu}
(e_i+1-\tau_i)(e_{i+1}+1-\tau_i)\\
\label{codBettib}&=\sum_{i>\mu}(\beta_i)(\nu_i).
\end{align}
\item\label{tauidealii} Assume that $k$ is algebraically closed.
 Then the closure
$\overline{G_\beta(H)}$ satisfies the frontier property,
\begin{equation}\label{frontierbetti}
\overline{G_\beta (H)}=\bigcup_{\beta'\ge\beta} G_{\beta'}(H).
\end{equation}
\item\label{tauidealiii} Assume that $k$ is algebraically closed. Then $\overline{\G_\beta(H)}$ is Cohen-Macaulay,
and has singular locus  $\bigcup_{\beta'>\beta} G_{\beta'}(H).$
\end{enumerate}
\end{theorem}
\begin{proof} 
 The irreducibility of each Betti stratum follows from the
Hilbert-Burch theorem \cite{Bu}. The frontier property will follow from our
characterization of the Betti strata, as an intersection of determinantal
varieties, and the closure properties of determinantal varieties. To
define the appropriate determinantal varieties we study the vector space $\langle
R_1I_i/yI_i\rangle$. 
\par
 By Proposition \ref{tauideala}, equation \eqref{eqbetamax}, Corollary \ref{betamax} and Lemma \ref{ranktau}, the
intersection of the rank one $(\tau_i=1)$ loci of $
\theta_\mu,\ldots ,\theta_{s-1}$ is the
$\beta_\max$ stratum $\G_{\beta_\max}(H)$ of $\G(H)$, and it has dimension
$\sum_{i>\mu}e_i =H_\mu$. Since $\dim \G(H)=\sum_{i\ge \mu} (e_i+1)(e_{i+1})$ by
 Theorem \ref{basic}\eqref{basici},
the stratum $\G_{\beta_\max}(H)$ has the expected codimension in $\G(H)$, namely 
\begin{equation}
\cod \G_{\beta_\max}(H)=\sum_{i\ge
\mu}^{s-1} (e_i)(e_{i+1}).
\end{equation}
 Thus, this locus is the
intersection of determinantal varieties, and the difference  on the left of the inequality
\eqref{codbstrata} in M.~Boij's Lemma \ref{boij3.2}
satisfies
\begin{equation*}
 \cod \G_{\beta_\max}(H) -\cod
X_{\beta_\max}(H)=0.
\end{equation*}
 This implies by \eqref{codbstrata} that for each $\beta$
the difference 
\begin{equation}\label{ecodi}
 \cod \G_{\beta}(H) -\cod
X_{\beta}(H)=0.
\end{equation} Recall that  here 
$\beta$ is defined by \eqref{etauequal} where $\tau$ satisfies
\eqref{etauinequal}, so $\G_\beta (H)\cap U_B$ is the subfamily parametrizing $I\in U_B$ satisfying
$\tau(I)=\tau$.
But by Lemma \ref{gencod} the codimension
$\cod X_\beta$ is just the expected codimension 
\begin{equation}
\cod X_\beta =\sum_{i>\mu ( H)}(\beta_i)(\nu_i)=\sum_{i\ge\mu (H) }(e_i+1-\tau_i)(e_{i+1}+1-\tau_i).
\end{equation}
It follows from \eqref{ecodi} that each intermediate stratum
$\G_\beta (H)\cap U_B$ has the expected codimension. This shows \eqref{codBettia}, and \eqref{codBettib} follows from 
\eqref{etauequal}. It follows from \eqref{tauilocus} that each
such locus $\G_\beta (H)\cap U_B$ is the product of 
the determinantal varieties given by the rank $\tau_i$ loci of $\theta_i$, as
claimed. \par When $k$ is algebraically closed, the frontier property 
and Cohen-Macaulayness of\linebreak $\overline{G_\beta (H)}\cap U_B$
and also the assertion in \eqref{tauidealiii} concerning the singularity locus follow from the analogous properties for
 a determinantal variety (see \cite{BV},\cite[\S 6.1]{W},\cite[II.4.1]{ACGH}), since $\overline{\G_\beta (H)}$ is a product of determinantal
 varieties.
 This
completes the proof of Theorem~\ref{tauideal}.
\end{proof}\par
The following corollary is the analogue of a result of S.~J.~Diesel  
 for height three Gorenstein algebras (see \cite[\S 3]{D} and also \cite[Theorem 5.25v.f]{IK}).
 We will denote by $[n]$ where $n$ is an integer, the
simply ordered set $\{0\le\ldots \le n\}$, with $n+1$ elements.  Let $H$ be an $O$-sequence. We denote by 
$\mathcal L(H)$ the following lattice that is the 
product of simply ordered sets,
\begin{equation}\label{eqlattice}
\mathcal L(H)=\prod_{\mu(H)\le i< s(H)}[e_i-(e_i-e_{i+1})^+].
\end{equation}
By \eqref{enorelsI} the simply ordered set $[e_i-(e_i-e_{i+1})^+]$ has cardinality the number of possible values for
$\beta_i$, given $H$.
We consider the map $\eta$ from the set of sequences $ \beta$ satisfying \eqref{enorelsI}
 to the set underlying $\mathcal L(H)$,
\begin{equation}
\eta(\beta)=(\beta_{\mu+1}-(e_\mu-e_{\mu+1})^+,\ldots ,\beta_{i+1}-(e_i-e_{i+1})^+,\ldots ,\beta_{s}-(e_{s-1}-e_{s})^+).
\end{equation}
We denote by $\beta (H)$ the lattice of sequences $\beta$ satisfying \eqref{enorelsI} under the termwise order of
Definition \ref{defbettistratum} --- this is the lattice of relation sequences $\beta$ possible for $H$.
\begin{corollary}\label{lattice} Let $k$ be an arbitrary field, and $H$ an $O$-sequence satisfying \eqref{Hcond}.
\begin{enumerate}[i.]
\item\label{maincori} The relation sequences
 $\beta$ that occur, for which
$\G_\beta(H)$ is nonempty are exactly those satisfying \eqref{enorelsI}. Likewise, each sequence $\nu$ satisfying \eqref{enogensI},
or sequence $\tau$ satisfying \eqref{etauinequal} occurs. 
\item\label{maincorii}The map $\eta$
 induces an isomorphism from the lattice $\beta(H)$ to $\mathcal L(H)$.
\item\label{maincoriii}
When $k$ is algebraically closed, the Zariski closure 
$\overline{G_\beta(H)}\supset G_{\beta'}(H)$ if and only if $\beta'\ge \beta$ in the lattice, and 
$\overline{G_\beta(H)}\cap G_{\beta'}(H)=\emptyset$
unless $\beta'\ge \beta$. 
\end{enumerate}
\end{corollary}
\begin{proof} When $k$ is finite, or $\cha k = p<s(H)$, the existence results \eqref{maincori} and \eqref{maincorii} are evident
from simple constructions of ideals having normal pattern; that these are the only sequences possible follows from 
Proposition \ref{tauideala}. The rest of the Corollary is immediate
from \eqref{enorelsI} and Theorem~\ref{tauideal}. 
\end{proof}
\begin{figure}[hbtp]
\begin{center}
\leavevmode 
\begin{picture}(35,45)(10,-5)
\setlength{\unitlength}{1mm}
\multiput(0,0)(0,7){2}{\line(1,0){35}}
\put(7,14){\line(1,0){21}}
\put(14,21){\line(1,0){14}}
\put(0,0){\line(0,1){7}}
\put(7,0){\line(0,1){14}}
\put(14,0){\line(0,1){21}}
\put(21,0){\line(0,1){21}}
\put(28,0){\line(0,1){21}}
\put(35,0){\line(0,1){7}}
\put(2.5,2.5){{\small \mbox{$1$}}}
\put(9.5,2.5){{\small \mbox{$x$}}}
\put(16,2.5){{\small \mbox{$x^2$}}}
\put(16,2.5){{\small \mbox{$x^2$}}}
\put(23,2.5){{\small \mbox{$x^3$}}}
\put(30,2.5){{\small \mbox{$x^4$}}}
\put(37,2.5){{\small \mbox{$x^5$}}}
\put(9.5,9.5){{\small \mbox{$y$}}}
\put(15.5,9.5){{\small \mbox{$yx$}}}
\put(22,9.5){{\small \mbox{$yx^2$}}}
\put(29.5,9.5){{\small \mbox{$yx^3$}}}
\put(16,16.5){{\small \mbox{$y^2$}}}
\put(22,16.5){{\small \mbox{$y^2x$}}}
\put(29.5,16.5){{\small \mbox{$y^2x^2$}}}
\put(22.5,23.5){{\small \mbox{$y^3$}}}
\end{picture}
\end{center}
\vspace{-.7cm}
\protect\caption{Normal pattern $B_H, H=(1,2,3,3,1,0)$, (see Example \ref{ex1}).}\label{fig1} 
\end{figure}
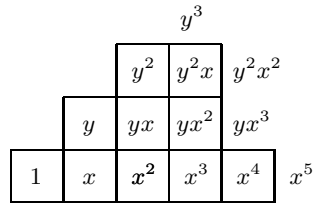
\begin{example}\label{ex1}Let $H=(1,2,3,3,1,0)$; then we have $\mu =3$, and the alignment character
\begin{equation*}
K=(3,3,1)^\ast +(2,1,0)=(3,2,2)+(2,1,0)=(5,3,2),
\end{equation*}
 the lengths of the rows of the normal pattern $B_H$ depicted
 in
Figure \ref{fig1}. The standard generators $F$ are (using
$a,b,c,d,e$ in place of $a_{3,2},\ldots $),
\begin{equation*}
f_0=x^5,\, f_1=yx^3+ex^4, f_2=y^2x^2+dx^4,
f_3=y^3+af_2:x+bf_1:x+cx^3.
\end{equation*}
The coefficients $a,b,c,d,e$ form the parameters $A$ for $U_B$. We have
$E(H)=(e_3,e_4,e_5)=(0,2,1)$, and $\dim \G(H)=(1\cdot 2+3\cdot 1)=5$, as expected. 
The map $\theta_4$ from $\langle xf_3,f_2,f_1\rangle $ to $\langle xf_1,f_0\rangle$, has matrix
\begin{equation*}
{N_4}=\left(
\begin{array}{ccc}
-d-ae&-e&1\\
c+ed-ad+ae^2&d-e^2&0
\end{array}
\right).
\end{equation*}
After natural
column operations as described in the proof of Theorem \ref{tauideal}, the matrix becomes
\begin{equation*}
{N_4}'=\left(
\begin{array}{ccc}
0&0&1\\
c+ed&d-e^2&0
\end{array}
\right).
\end{equation*}
 The map $\theta_3$ has matrix ${N_3}={N_3}' = (1,0,0)^T$, of rank one.
The rank one locus of $\theta_4$ is given by $d=e^2,c=-ed=-e^3$, so has codimension two,
and dimension three equal to $H_\nu=H_3$ (see Corollary \ref{betamax}). Also,
by equation \ref{codBettia} of Theorem \ref{tauideal} for $\tau_4=1$, this locus has codimension $(3-1)\times (2-1)=\beta_4\times \nu_4=2$.
The only other stratum is the rank two locus of $\theta_4$, which is the generic stratum of dimension 
five, open dense in
$\G(H)$.
\end{example}
\begin{example}\label{ex2} Let $H=(1,2,3,4,2,1,0)$. Then $E(H)=(e_4,e_5,e_6)=(2,1,1)$, and $\G(H)$ has dimension
$(3\cdot 1+2\cdot 1)=5.$ The standard generators $F$ are (using
$a,b,c,d,e$ in place of $a_{4,2},\ldots $),
\begin{equation*}
\begin{split}
f_0&=x^6,\, f_1=yx^4+ex^5, f_2=y^2x^2+cf_1:x+dx^4,\\
f_3&=y^3x+b f_1:x+\alpha_1x^4, f_4=y^4+a f_1:x+\alpha_2 x^4,
\end{split}
\end{equation*}
where $a,b,c,d,e$ are the parameters $A$, and $\alpha_1,\alpha_2$ may be written as 
polynomials in the
parameters. The non-trivial $\theta$ maps are $\theta_4$ and $\theta_5$. The matrix $N'_4$ for the map $\theta_4$ from $\langle f_4,f_3,f_2\rangle\to
\langle xf_2,f_1\rangle$ satisfies
\begin{equation*}
{N_4}'=\left(
\begin{array}{ccc}
0&0&1\\
a&b&0
\end{array}
\right) \mod A^2.
\end{equation*}
The matrix $N'_5$ for the map $\theta_5$ from 
$\langle xf_2,f_1\rangle\to \langle xf_1,f_0\rangle$ satisfies
\begin{equation*}
{N_5}'=\left(
\begin{array}{cc}
0&1\\
d-ce+e^2&0
\end{array}
\right).
\end{equation*}
The Betti strata of $\G(H)$ are determined by $\tau(I_4)=$ rank $(\theta_4)$, and  $\tau (I_5)$ which 
 can each be one or two. So there are four Betti strata of $\G(H)$, of codimensions
zero for the stratum ($\G_{\beta(0)}(H): \ \tau_4=\tau_5=2$), one for ($\G_{\beta(1)}(H)\ \tau_4=2,\tau_5=1$),
 two for ($\G_{\beta(2)}(H):\ \tau_4=1,\tau_5=2$), and three for the stratum
($\G_{\beta_\max}(H):\ \tau_4=\tau_5=1$). The values of $(\beta_{5},\beta_6)$
that correspond are $\beta(0): (1,0)$, $\beta(1): (1,1)$, $\beta(2): (2,0)$, $\beta_\max : (2,1)$.
By Corollary \ref{lattice} the lattice $\mathcal L(H)=[1]\times [1]$, so is not simply ordered.
In fact, $\overline{G_{\beta(1)}(H)}=G_{\beta(1)}(H)\cup G_{\beta_\max}(H)$ and
 $\overline{G_{\beta(2)}(H)}=G_{\beta(2)}(H)\cup G_{\beta_\max}(H)$.
\end{example}
\section{Algebra quotients of $k[x,y]$ having fixed socle type}\label{fixed}
In this section we characterize the Hilbert functions consistent with a given socle type (Theorem
\ref{4.6C}), and we also determine the Hilbert function of the intersection of a general
enough set of $t$ level ideals, each having a specified Hilbert function 
(Theorem \ref{intersectlevelideals}). \par.
\begin{definition}\label{socletype} The \emph{socle type} $\ST$ of an
Artinian algebra
$A=R/I$ is the Hilbert function of the \emph{socle} $\Soc
(A)=\langle 0:M\rangle\subset A$. 
\end{definition}
 The following result
was announced as
\cite[Theorem 4.6C]{I3}. The case $\ST$ concentrated in a single degree
$j$ (so $\ST _i=0
$ unless
$i=j$), announced in
\cite[Theorem 4.6A]{I3}, was shown in \cite[Theorems 2.17,2.19]{I4} for
$\LA_N(d,j)$ (see Definition \ref{deflevel}), where the resulting Hilbert function was
related to a partition $P$.
\begin{theorem}\label{4.6C} Let the sequence 
$H$ satisfy \eqref{Hcond}, and $\lim_{i\to\infty}H_i=0$, and suppose $k$ is arbitrary. 
The sequence
$H$ occurs for the Hilbert
function of a graded algebra of socle type $\ST$ if and only if for all
$i$, 
\begin{equation}\label{WandH}
e_{i+1}-e_{i+2}\le \ST _i\le e_{i+1}.
\end{equation}
\end{theorem}
\begin{proof}  By
\eqref{Socandrel}, that
$H$ and
$\ST$ must satisfy
\eqref{WandH} is equivalent to
 \begin{equation*}
\beta_{i+2}\le e_{i+1}\le
e_{i+2}+\beta_{i+2},
\end{equation*}
 which follows from
\eqref{enorelsI} and \eqref{Hgenrel}. Suppose conversely 
$\ST$ satisfies
\eqref{WandH}. Setting
$\beta_i=\ST_{i-2}$ and letting
$\nu_{i+2}=\ST_i-(e_{i+1}-e_{i+2})$, we have sequences $\beta$ and $\nu$
satisfying the conditions \eqref{enorelsI} or \eqref{enogensI}, so by
Corollary \ref{lattice}\eqref{maincori} they occur.
\end{proof}
\begin{lemma}\label{missA} Let $H$ satisfy the conditions of
\eqref{Hcond}, with $\lim_{i\to\infty}H_i=0$,  and let
$W=\oplus W_i, W_i\subset R_i$ where $
\mu(H)\le i\le s-1$ be given vector spaces of dimension $\dim_k
W_i\le H_i$. Then there is an ideal
$I$ satisfying both $I_i\cap W_i=0 $ and $\tau (I_i)=1\,\forall i \mid
\mu\le i\le s-1$. Given a  sequence $\beta$ satisfying
\eqref{enorelsI}, there is an open dense subset of ideals
$I\in G_\beta (H)$ satisfying $I\cap W=0$.
\end{lemma}
\begin{proof} The last assertion is a consequence of the first assertion and 
Theorem \ref{tauideal}\eqref{tauidealii}, as the stratum $H_{\beta_\max}$ is in the closure of 
$H_\beta$, and the condition that $I$ miss $W$ is an open condition on Artinian algebras
$A=R/I$ in $H_\beta$. For the first assertion we may assume
without loss of generality that
$\dim_k W_i = H_i$. Given the characterization of the extremal Betti stratum
in the proof of Proposition~\ref{tauideala}, it suffices to show that given $W_i$ of dimension $c=
H_i$, there is a vector space of the form $fR_{i-c}, \deg f=c$, missing
$W_i$. Dualizing, we have
\begin{equation*}
W_i\cap fR_{i-c}=0\Leftrightarrow W_i^\perp+\langle
fR_{i-c}\rangle^\perp=\mathcal R_i.
\end{equation*}
Let $f=\prod_{v=1}^c (x-a_vy)$, 
with distinct $a_v$. Then we have (see, for example, \cite{EmI} or
\cite[Lemma 1.15]{IK})
\begin{equation}\label{eannihilators}
\langle fR_{i-c}\rangle^\perp =\langle L_1^i,\ldots L_c^i\rangle, \text {
where } L_v=a_vX+Y.
\end{equation}
 It is well known and classical that one can assure that
$W_i^\perp \cap L=0$, by choosing for $L$ a complementary-dimension vector
space spanned by
$i$-th powers of linear forms (this is used extensively in \cite{I3}).
Thus, one may build $I$ as follows, beginning in degree
$i=(s-1)$ taking $f_s=1$, and continuing down to degree $i=\mu$. We choose
as the first part of the $i$-labelled step,
$e_{i+1}$ linear forms such that for each degree $u\le i$ the $u$-th
powers of the $H_i=e_{i+1}+\cdots +e_s$ forms chosen so far are linearly
independent from
${W_u}^\perp$. In the second part of each step, we determine $f_i$ by
multiplying the previously chosen $f_{i+1}$ by the
degree $e_{i+1}$ homogeneous form whose linear factors are the
annihilators as in \eqref{eannihilators} of the $e_{i+1}$ linear forms
chosen in the first part of the step. 
\end{proof}\par
\begin{definition}\label{deflevel}
An algebra $A=R/I$ is a \emph{level algebra} if its socle
 $\Soc (A)=0:M$ is concentrated in a single degree $j$. The level algebra  $A(V)$ determined by
a vector subspace $V\subset R_j$ is the quotient $A(V)=R/L(V)$ of $R$ by the \emph{level ideal} $ L(V)=\overline{V}+M^{j+1}$, satisfying
\begin{equation}L(V)=R_{-j}V\oplus\cdots \oplus R_{-1}V\oplus V\oplus M^{j+1}.
\end{equation}
We denote by $\LA_N(d,j)$ the family of level algebras having socle type $d$
 in degree $j$, and Hilbert function $N$. See \cite{I4} and for a geometric approach \cite{ChG}.
\end{definition}
It follows from \eqref{Socandrel} that $\LA_N(d,j)$ is just the open dense subscheme Betti stratum of $\G(N)$, parametrizing
those ideals having the minimum possible
Betti numbers.
We recall the following characterization of level Hilbert functions $N$ from
\cite[Lemma 2.15]{I4}:
\begin{lemma}\label{NT} The Hilbert function $N=(1,\ldots ,N_j,0)$ of a level algebra
$\LA (V)$ determined by the vector subspace $ V\subset R_j, \dim_k V =d $
satisfies 
\begin{align}\label{eN}
 N_j&=j+1-d, N_i=0 \text { for } i>j, \text { and }\notag\\
e_{j+1}(N)&=j+1-d\ge e_j(N)=\tau (V)-1\ge e_{j-1}(N)\ge \cdots ,
\end{align}
where $1\le \tau (V)\le \min\{
d,j+2-d\}$.
\end{lemma}
The following result characterizes the Hilbert function of the intersection of general enough level algebras. The 
 special case for \emph{compressed algebras}
--- when each $N(i)$ is the maximum Hilbert function given $(d_i,j_i)$ --- has an analogue
in embedding dimension $r$ that is a consequence of \cite[Proposition 3.6]{I3}: 
for the statement, one replaces $u+1$ in \eqref{hilbintlevel} by $\dim_k R_u$. Likewise, 
there is an analogue in embedding dimension $r$ for the slightly more general case of \emph{power sum algebras}, 
whose Hilbert functions $N(i)$ are the maximum bounded by $t_i$, given $d_i,j_i$ 
\cite[Theorem 4.8B]{I3}. What is unusual in our result here is that there is no restriction on the level
algebra Hilbert functions $N(i)$. Usually, it is difficult to build algebras of expected
``sum'' Hilbert functions from intersecting a set of ideals of 
given Hilbert functions: it is difficult to show that
 the intersections are dimensionally proper. Of course, it is the embedding dimension
two that helps us.
\begin{theorem}\label{intersectlevelideals} Let $R=k[x,y]$. Fix Hilbert functions
$N(1),\ldots ,N(t)$, and sequences of integers
$j_1\ge \ldots\ge j_t$, $d_1,\ldots d_t$ where $1\le d_i\le j_i+1$, such
that each
$N(i)$ satisfies \eqref{eN} of Lemma \ref{NT} for the Hilbert function of
level algebras $L(V(i))$
determined by a vector subspace $V(i)$ of dimension $d_i$ and degree
$j_i$ in $R=k[x,y]$. Let
$I(1),I(2),\ldots ,I(t)$, where
$I_i=L(V_i)=\overline{V_i}+M^{j_i+1}, V_i\subset 
R_{j_i}, j_1\ge j_2\ge
\cdots \ge j_t$ be \emph{general enough} level ideals of
Hilbert functions $H(R/I(i))=N(i)$. Then the intersection $I=I(1)\cap
\cdots \cap I(t)$ satisfies
\begin{equation}\label{hilbintlevel}
H(R/I)_u=\min\{ \sum H(R/I(i))_u,u+1\}.
\end{equation}
 For general enough
$I(1),\ldots ,I(t)$ the $\tau$ sequence $\tau(I)$ for the intersection
$I$ satisfies for
$u\ge
\mu (I)$ 
\begin{equation}
\tau (I_u)=\sum_{i\ge u} \tau(V_i).
\end{equation}
\end{theorem}
\begin{proof}The key issue is to satisfy \eqref{hilbintlevel}, which is
the condition that the vector spaces $I(1)_u,\ldots ,I(t)_u$ intersect
properly in $R_u$ for each $u$. This is an open condition on the product
$P=\prod_i \LA_{N(i)}(d_i,j_i)$, an irreducible variety by \cite[Theorem 2.17]{I4}. So it suffices to produce a single algebra having that
Hilbert function and the given socle type, and such that successively, the
algebras determined by $N(1)$, by $\{ N(1),N(2)\},\ldots $ (beginning in
the highest socle degree)  have the expected Hilbert functions.
We accomplish this by successive applications of Theorem
\ref{tauideal} and Lemma \ref{missA}.
\end{proof}
\begin{remark}{ \sc Betti strata for codimension two punctual schemes
in
$\mathbb P^2$, and ACM curves in $\mathbb P^3$.}\label{higheremb} The codimension formula of
Theorem
\ref{tauideal} is also valid for the Betti strata of the postulation stratum
$\Hilb^T(\mathbb P^2)$  of the punctual Hilbert scheme
$\Hilb^n(\mathbb P^2),n=\sum H_i, H=\Delta T$ (here $H$ is termed the $h$-vector of $T$), since a punctual scheme is
arithmetically Cohen-Macaulay (ACM). The proof can use Gotzmann's result
that $\Hilb^T (\mathbb P^2)$ is fibred over $\G(H)$ by an affine space of
fixed dimension
\cite{Go} (see \cite[\S 5.5, especially p. 180]{IK}).
A well known result of L. Gruson and C.~Peskine, and R. Maggioni and
A.~Ragusa, is that a reduced, irreducible ACM curve
$X\subset
\mathbb P^3$ has a \emph{numerical character} without gaps, or,
equivalently, has $H$-vector  $H_X=\Delta^2 H(k[x_0,x_1,x_2,x_3]/I_X)$
that is \emph{decreasing}: for $i\ge \mu(H), t_i<\mu\Rightarrow
t_i>t_{i+1}$: that is the $h$-vector
$H$ is strictly decreasing after its first decrease from its maximum
value of $\mu$ \cite{GrP,MaR}. \par
Considering ACM curves in $\mathbb P^3$, T. Sauer has shown that each Betti stratum
possible for an Artinian algebra $A=k[x,y]/I$ having decreasing Hilbert
function $H$ occurs for a smooth ACM curve X in $\mathbb P^3$ \cite{Sa}.
J. Herzog, N.~V.~Trung and G. Valla show further that a reduced ACM
variety $Y\subset \mathbb P^n$ of codimension two is a hyperplane section
of a reduced irreducible CM normal variety $X\subset \mathbb P^{n+1}$ of
codimension two if and only if the degree-matrix of the minimal resolution satisfies
a condition equivalent to the $H$-vector being decreasing \cite[Theorem
1.1]{HTV}. See
\cite{GM} for a careful discussion of these results. It follows from
the above mentioned results of T. Sauer et al that the codimension formula of Theorem
\ref{tauideal} for the special case $H$ decreasing will also give the
codimension of the corresponding Betti substratum of the postulation
scheme
$\Hilb_{\mathrm{Sm}}^{d,g,T}(\mathbb P^3)$ parametrizing smooth ACM
curves in $\mathbb P^3$, of postulation determined by the $h$-vector $H=\Delta T$,
degree
$d=\sum H_i$ and genus
$g=\sum_{i\ge 2}(i-1)H_i$. The postulation stratum for smooth ACM
curves, as well as its Betti substrata are irreducible. Of
course, there can be other irreducible components of
$\Hilb^{d,g}(\mathbb P^3)$ corresponding to non-ACM curves  \cite{PS}.
\end{remark}
\begin{ack}
We acknowledge gratefully discussions with Mats Boij, whose method in
studying the Betti strata of height three Gorenstein algebras we apply
here. Discussions with Jacques Emsalem and Vassil Kanev were helpful
 background for this article. We are grateful also for the referee's comments, and to a referee for \cite{I4} who
had suggested making this a separate article.
\end{ack}
\bibliographystyle{amsalpha}

\end{document}